# The attractor structure of logarithmic iterations in the complex plane

Pascal Wallisch
Center for Neural Science
New York University, New York, NY

**Abstract**
*We use the methods of empirical mathematics to show that iterative logarithmic operations will result in an attractor point on the complex plane. Moreover, we demonstrate that different bases converge onto different attractors. Finally, we elicit the geometric structure of these attractors*

**Introduction**
This paper represents an exercise in empirical mathematics. In other words, we try to elicit the theoretical properties of a mathematical construct by observation, while varying properties that are believed to systematically affect the phenomena in question, much like one would study an empirical phenomenon in the outside world. The crucial difference to actual empirical science is the space in which the "measurements" occur. In empirical science, the data-generating phenomena in question inhabit the external world, or rather physical, biological or chemical aspects of it. By virtue of measurement theory, these data are then brought into the realm of mathematics[1,2,3]. Mathematical analysis and simulation are then used to treat these mathematical entities (usually numbers), the conclusions of which are then used to inform our understanding of the empirical phenomenon in the outside world. By this method, empirical science gains knowledge about the truth status of synthetic propositions – those that are not necessarily true from their antecedents. Conversely, pure mathematics is concerned with the truth status of analytic propositions – those that are necessarily true, given their axioms, the implications of which have to be worked out through logic[4]. Both realms of inquiry have a long tradition, strong contemporary followings and can boast a rich hoard of impressive successes. However, there is another – oft neglected – possible form of inquiry, namely empirical mathematics. As pointed out above, empirical mathematics uses the methods of empirical science to elicit the truth status of synthetic propositions. The crucial difference being that the truth statements concern phenomena not in the outside world, but in the realm of mathematics itself[5,6]. These statements are synthethic insofar as their truth status is not deduced from the axioms, but rather from an observation of the mathematical "data" in the style of empirical science. In a sense, this approach is more elegant than that of genuine empirical science, as we can forego the need for measurement theory and operationalization, both of which are problematic insofar as they are not unique, in ontological terms[7]. This approach is not without precedent. Historically, this was practiced for instance in Babylonian geometry[8,9], but has since – due to the dramatic success of greek-style mathematics based on axioms, proofs and deduction – largely foundered. It is now time to reconsider the approach – on the one hand, many mathematical phenomena of contemporary interest are too complex to attempt closed form solutions. This impression is closely related to the notion that many phenomena of interest to empiricists do not lend themselves to analytical solutions any more – instead, a numerical or computational solution is sought[10]. On the other hand, computational power continues to grow with no sign of abating. This suggests that the time might be ripe for a reconsideration of this – long neglected – approach. The approach is also implicit in other recent attempts to enrich the scientific landscape[11].
In theoretical terms, our claim, hypothesis or synthetic proposition is very simple: If one iteratively applies logarithmic operations, the result will converge onto an attractor point on the complex plane, regardless of starting position. Moreover, logarithms with different bases converge onto different attractors on the plane in a very systematic manner, forming a circle in this plane. We will use empirical mathematics to test this assertion.

## Methods and Results

All calculations were performed with Matlab® (Mathworks Inc.). The basic idea is illustrated in figure 1. Iteratively taking the natural logarithm from an arbitrary starting point in the real numbers (say the number 5) leads to an attractor at 0.3181 + 1.3372i. In this figure, as in all subsequent figures, the abscissa represents the real part while the ordinate represents the imaginary part of the complex number. Note how the attractor is approached in the manner of a logarithmic spiral. Table 1 represents the data points that correspond to those in figure 1, illustrating the iterations from 1 to 50 and their corresponding numeric values.

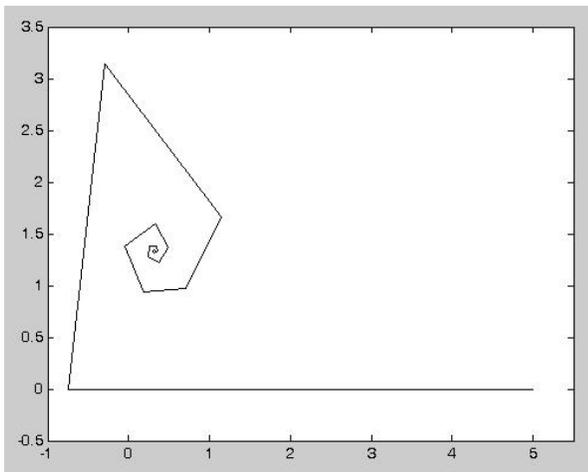

Figure 1

Table 1

| | | | | |
|---|---|---|---|---|
| 1 | 5 | 26 | 0.319 + | 1.336i |
| 2 | 1.6094 | 27 | 0.3174 + | 1.3364i |
| 3 | 0.47588 | 28 | 0.31743 + | 1.3376i |
| 4 | -0.74258 | 29 | 0.31828 + | 1.3378i |
| 5 | -0.29763 + | 3.1416i | 30 | 0.31856 + | 1.3372i |
| 6 | 1.1492 + | 1.6653i | 31 | 0.3182 + | 1.3369i |
| 7 | 0.70473 + | 0.96674i | 32 | 0.31793 + | 1.3371i |
| 8 | 0.17927 + | 0.94089i | 33 | 0.31803 + | 1.3374i |
| 9 | -0.043099 + | 1.3825i | 34 | 0.3182 + | 1.3373i |
| 10 | 0.3244 + | 1.602i | 35 | 0.31821 + | 1.3372i |
| 11 | 0.49132 + | 1.371i | 36 | 0.31812 + | 1.3372i |
| 12 | 0.37595 + | 1.2267i | 37 | 0.31809 + | 1.3372i |
| 13 | 0.2492 + | 1.2734i | 38 | 0.31812 + | 1.3373i |
| 14 | 0.26049 + | 1.3775i | 39 | 0.31815 + | 1.3372i |
| 15 | 0.33787 + | 1.3839i | 40 | 0.31814 + | 1.3372i |
| 16 | 0.35386 + | 1.3313i | 41 | 0.31812 + | 1.3372i |
| 17 | 0.32032 + | 1.311i | 42 | 0.31812 + | 1.3372i |
| 18 | 0.29979 + | 1.3312i | 43 | 0.31813 + | 1.3372i |
| 19 | 0.31079 + | 1.3493i | 44 | 0.31814 + | 1.3372i |
| 20 | 0.32542 + | 1.3444i | 45 | 0.31813 + | 1.3372i |
| 21 | 0.32442 + | 1.3333i | 46 | 0.31813 + | 1.3372i |
| 22 | 0.31642 + | 1.3321i | 47 | 0.31813 + | 1.3372i |
| 23 | 0.31421 + | 1.3376i | 48 | 0.31813 + | 1.3372i |
| 24 | 0.31772 + | 1.3401i | 49 | 0.31813 + | 1.3372i |
| 25 | 0.32007 + | 1.338i | 50 | 0.31813 + | 1.3372i |

Interestingly, this process converges on the same point after sufficient iterations, regardless of the starting point – for a given base. This is demonstrated in figure 2 for start values of [-5 -4 -3 -2 -1 0 1 2 3 4 5] in the real part fully crossed with [1 2 3] in the complex part. Figure 2 a shows the full figure and figure 2 b shows a magnification of the part around the attractor point.

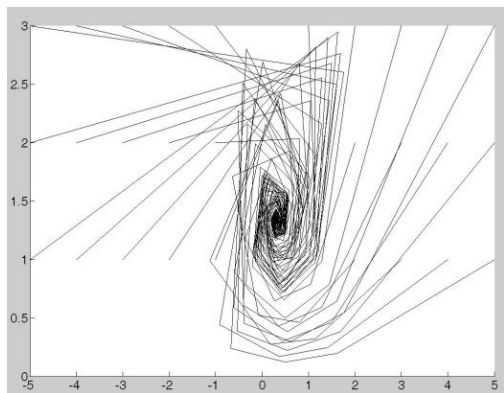
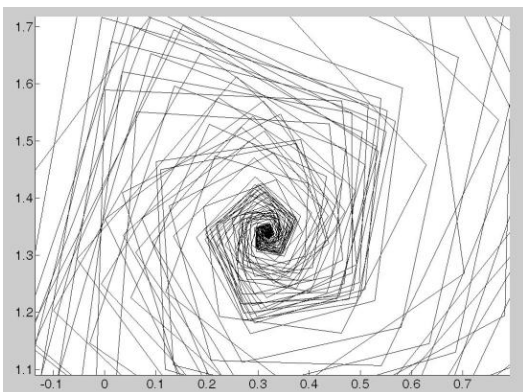

Figure 2a

Figure 2b

Note that the actual iteration values are connected with a line for clarity in figures 1 and 2. However, not all logs converge to a value around 0.31813 + 1.3372i. Each base has a distinctly unique attractor, as exemplified in figure 3. Figure 3a exemplary shows the attractors of bases 2 (green), e (blue) and 10 (red). In this case, the attractor is approached from different starting points in the reals, specifically all numbers from 0 to 10, in increments of 0.1. In this case, the points are not connected, but represent the superimposition of 101 iterations for each log base. Figure 3b adds points for base 1.6 (black)

Figure 3a

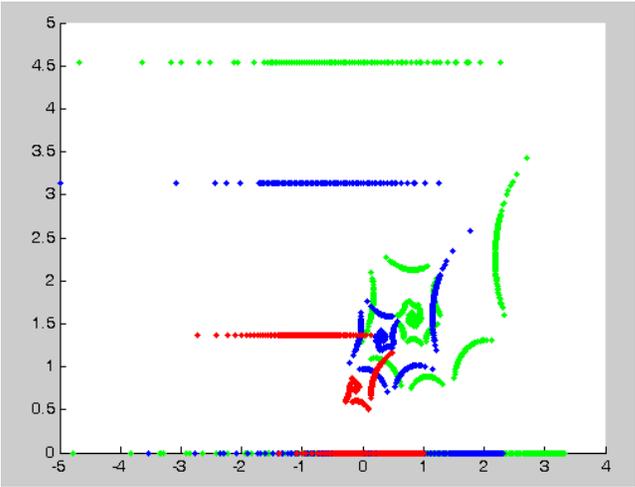

Figure 3b

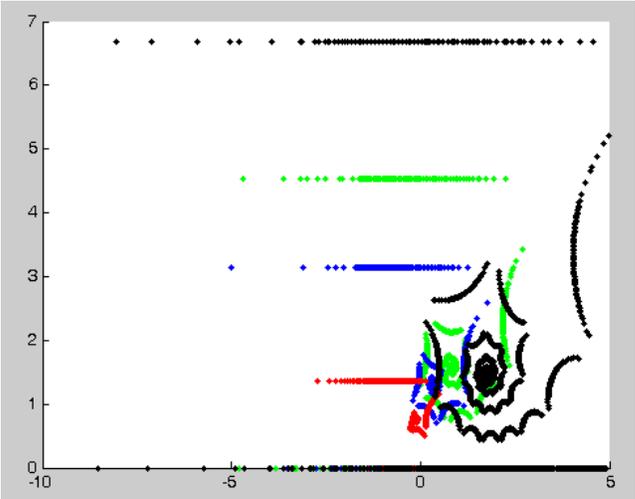

From the figures above, it should be clear that the attractor point is systematically dependent on the log base, but not the starting point of the iteration sequence.
Table 2 illustrates the attractor points of the log bases illustrated in Figure 3b.

| Log base | Attractor point |
|---|---|
| 10 | -0.1192 + 0.7506i |
| E | 0.3181 + 1.3372i |
| 2 | 0.8247 + 1.5674i |
| 1.6 | 1.7792 + 1.4699i |

Table 2

This progression immediately begs the question if these attractor points form a structure in the complex plane as a function of the log base.

As it turns out, this is in fact the case. Figure 4 shows attractor points in the complex plane (the abscissa is the real part, the ordinate is the imaginary part) for all real bases from positive infinity (which yields attractor points close to 0) to 1.38, which yields an attractor point near the real number 5. If the bases continue towards 1, the attractor point shifts to large real positive numbers onwards to positive infinity. In between these extremes, the attractor points form a half-circle in the complex plane. Note that the chosen aspect ratio makes the half- circle look elliptical in figure 4.

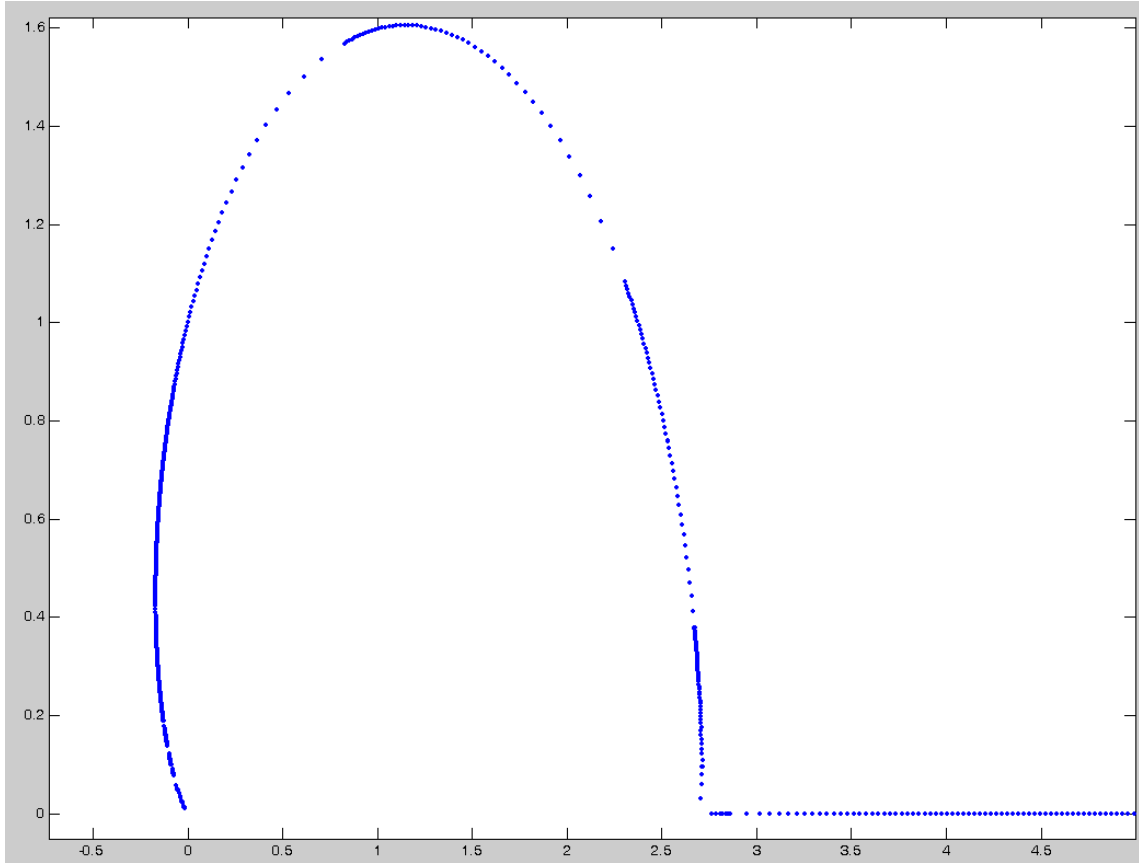

Figure 4

See table 3 for some crucial attractor values in Figure 4, relative to the corresponding base.

| Log base | Attractor point | Log base | Attractor point |
| --- | --- | --- | --- |
| 1.01 | 651.1 | 1.8 | 1.164 + 1.6045i |
| 1.1 | 38.23 | 2 | 0.8247 + 1.567i |
| 1.2 | 14.77 | e | 0.3181 + 1.3372i |
| 1.3 | 7.86 | 5 | -0.011 + 0.983i |
| 1.3795 | 5 | 10 | -0.12 + 0.75i |
| 1.4 | 4.41 | 20 | -0.155 + 0.61i |
| 1.44 | 3.12 | 50 | -0.169 + 0.487i |
| 1.444 | 2.86 | 100 | -0.17 + 0.424i |
| 1.44465 | 2.74 | 1000 | -0.157 + 0.2978i |
| 1.444668 | 2.7178 | 1e6 | -0.117 + 0.1596i |
| 1.444669 | 2.7275 + 0.00005i | 1e9 | -0.09 + 0.1098i |
| 1.45 | 2.67 + 0.38i | 1e12 | -0.07 + 0.084i |
| 1.475 | 2.48 + 0.85i | 1e25 | -0.0478 + 0.042i |
| 1.5 | 2.31 + 1.08i | 1e50 | -0.02885 + 0.022i |
| 1.75 | 1.28 + 1.60i | 1e100 | -0.0169 + 0.0111i |

Table 3

Note the sudden liftoff in the imaginary parts at an attractor point value of around real e. The maximum i value is reached somewhere around a log base of 1.8.

Naturally, this invites the question if this is the entire structure of attractor points in the complex plane. A half-circle is suggestive. Can we complete it towards a full circle? Thus far, all starting points that we examined yielded attractor points in the reals or positive complex numbers. The answer to this question is simple: As long as one starts in the reals of positive complex numbers, one will stay in them, given a log base between 1 and positive infinity. However, if one starts in the negative complex numbers (with negative i), the attractor points form a mirror image in the negative complex numbers. This is illustrated in figure 5, where we again show the attractor point for starting point 5 + i and 5 − i, both with a log base of e. For effect, the points between starting point and attractor point are connected with lines. As usual, abscissa is real parts whereas ordinate is imaginary parts of the complex number.

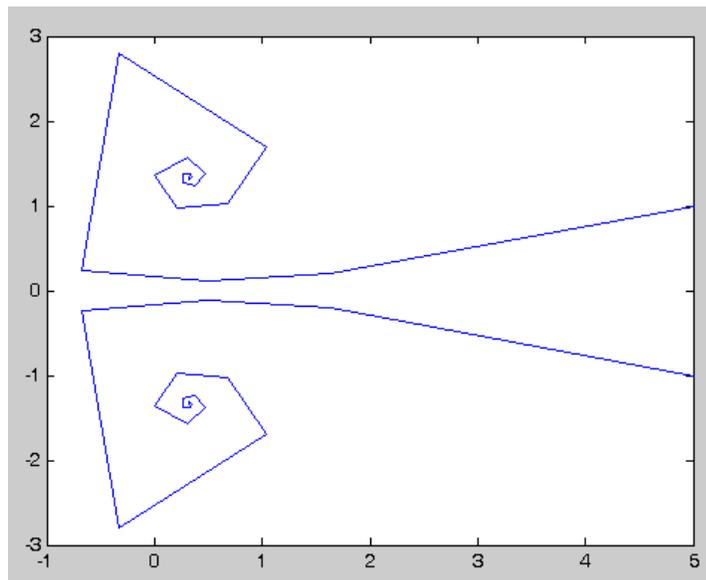

Figure 5

As pointed out above, considering this additional domain of starting points yields a mirror image of the attractor points for positive complex numbered starting points. In effect, this both completes the attractor half "circle" discussed above and reveals a bifurcation point at log bases of slightly larger than real 1.444668. This is illustrated in figure 6 with the points from table 3. Starting points are +i and –i, respectively.

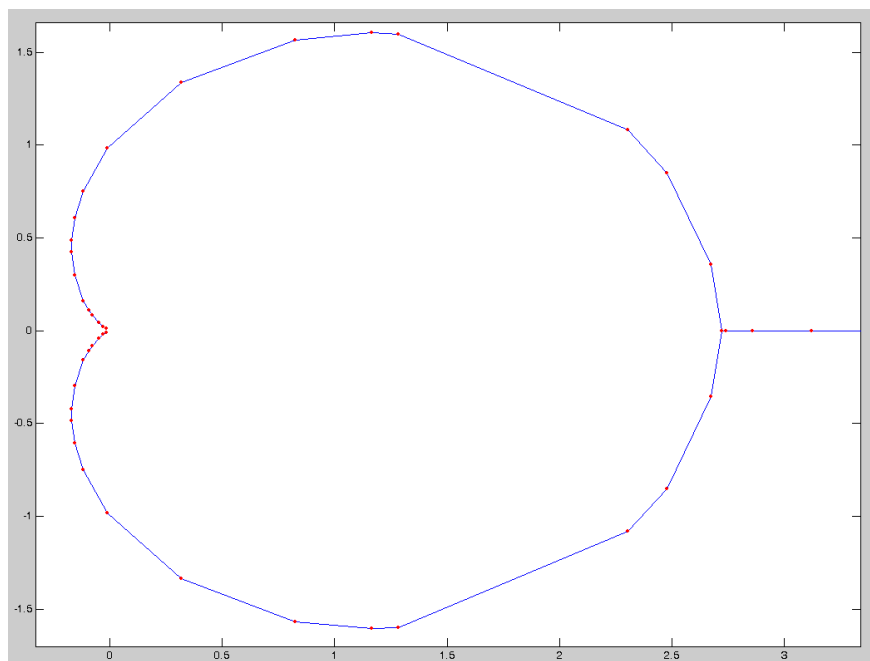

Figure 6

Of course, the figure could be improved by finer sampling, particularly between 1.4 and 2. This is illustrated in figure 7.

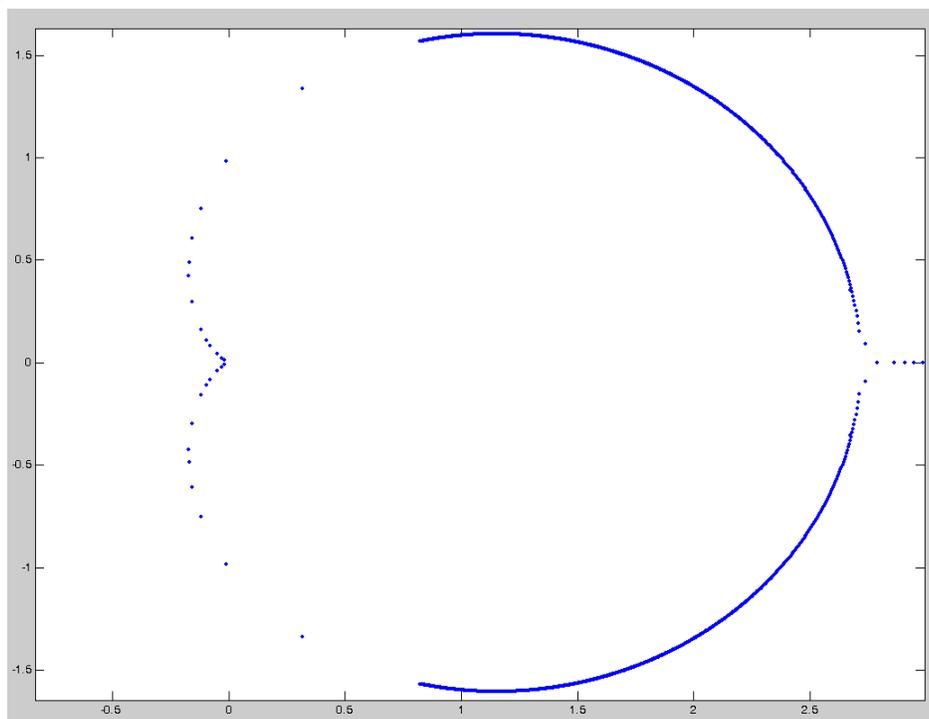

Figure 7

This leaves the question what happens in cases of log bases smaller than 1. Behavior like that is illustrated in figure 8 for log bases between 0.999 and 0.01 in steps of 0.001. It turns out, the attractor in this case is approximately linear. The attractors fall along a line towards 0, given starting points in positive i. The attractors fall along a line rising towards 0, given starting points in negative i.

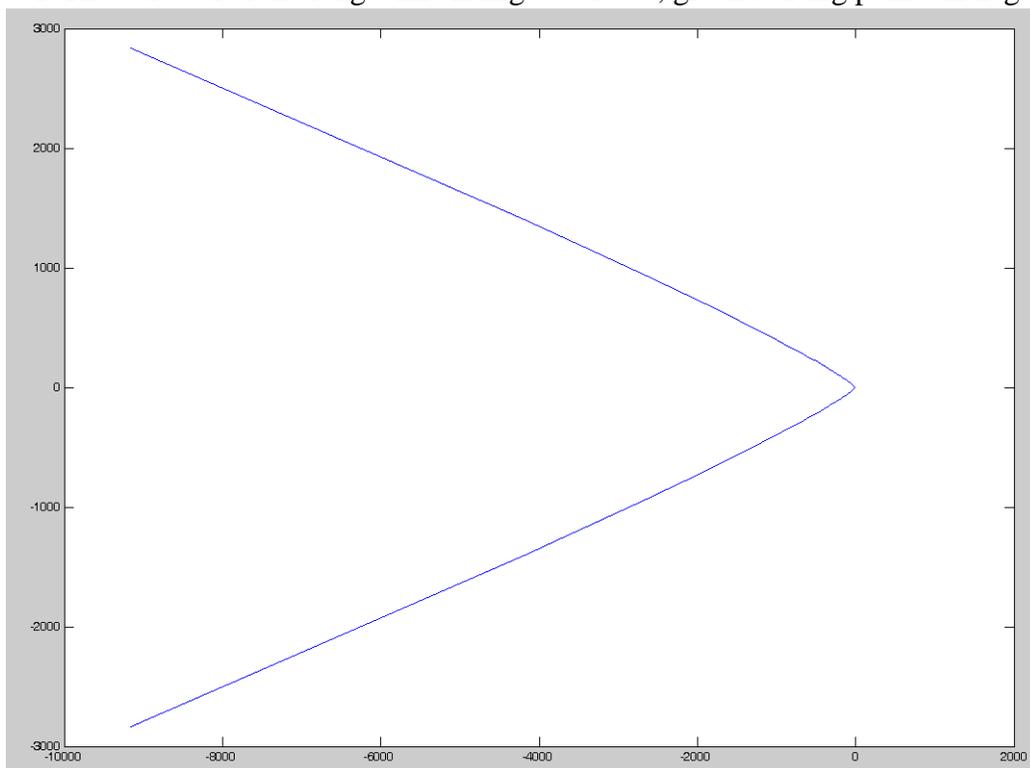

Figure 8

What about attractors for negative log bases?

Incrementing log bases from -0.001 to -2 in steps of 0.001 yields figure 9. Interestingly, both start domains (positive and negative i) collapse into the same curve. As one steps into the negative bases, the attractor curve bends towards positive i's.

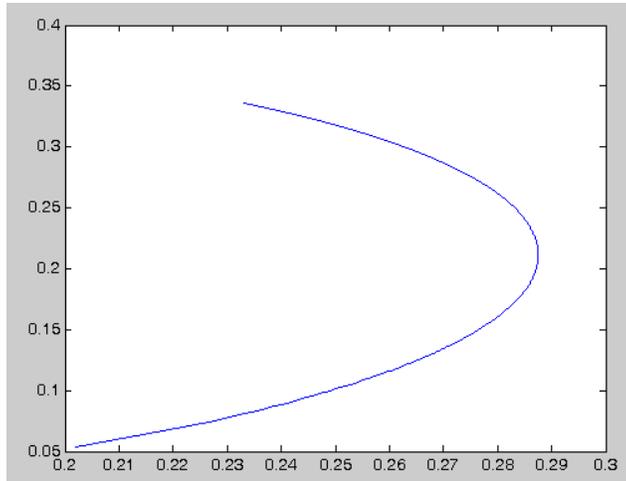

Figure 9

Finally, one can raise the issue of imaginary or complex bases. Imaginary bases indeed also have a distinct structure. For example, we just show 10 points (1i to 10i) in figure 10. These points are illustrated in table 4.

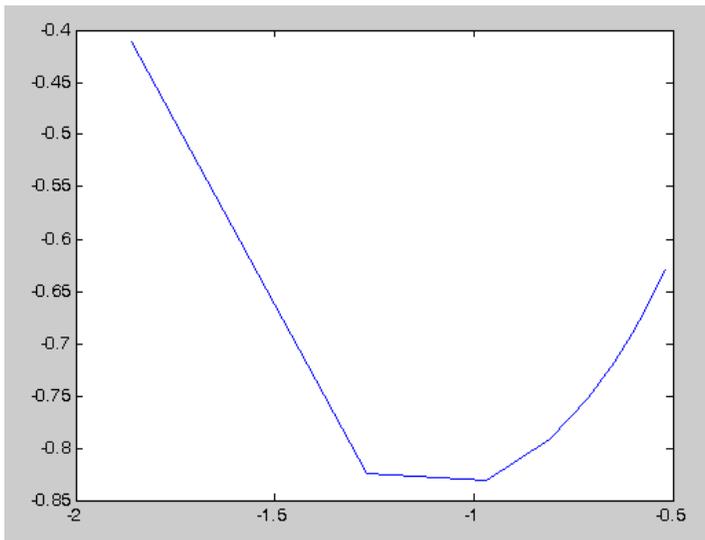

| Base | Attractor point |
|------|-----------------|
| I    | -1.8617 - 0.4108i |
| 2i   | -1.2697 - 0.8243i |
| 3i   | -0.9670 - 0.8309i |
| 4i   | -0.8080 - 0.7915i |
| 5i   | -0.7122 - 0.7520i |
| 6i   | -0.6481 - 0.7182i |
| 7i   | -0.6022 - 0.6899i |
| 8i   | -0.5674 - 0.6661i |
| 9i   | -0.5400 - 0.6458i |
| 10i  | -0.5179 - 0.6283i |

Table 4

Figure 10

All of these structures elicited can be combined in one figure, figure 11. Abscissa and ordinate indicate the real and imaginary part of the complex number that corresponds to the attractor, respectively.

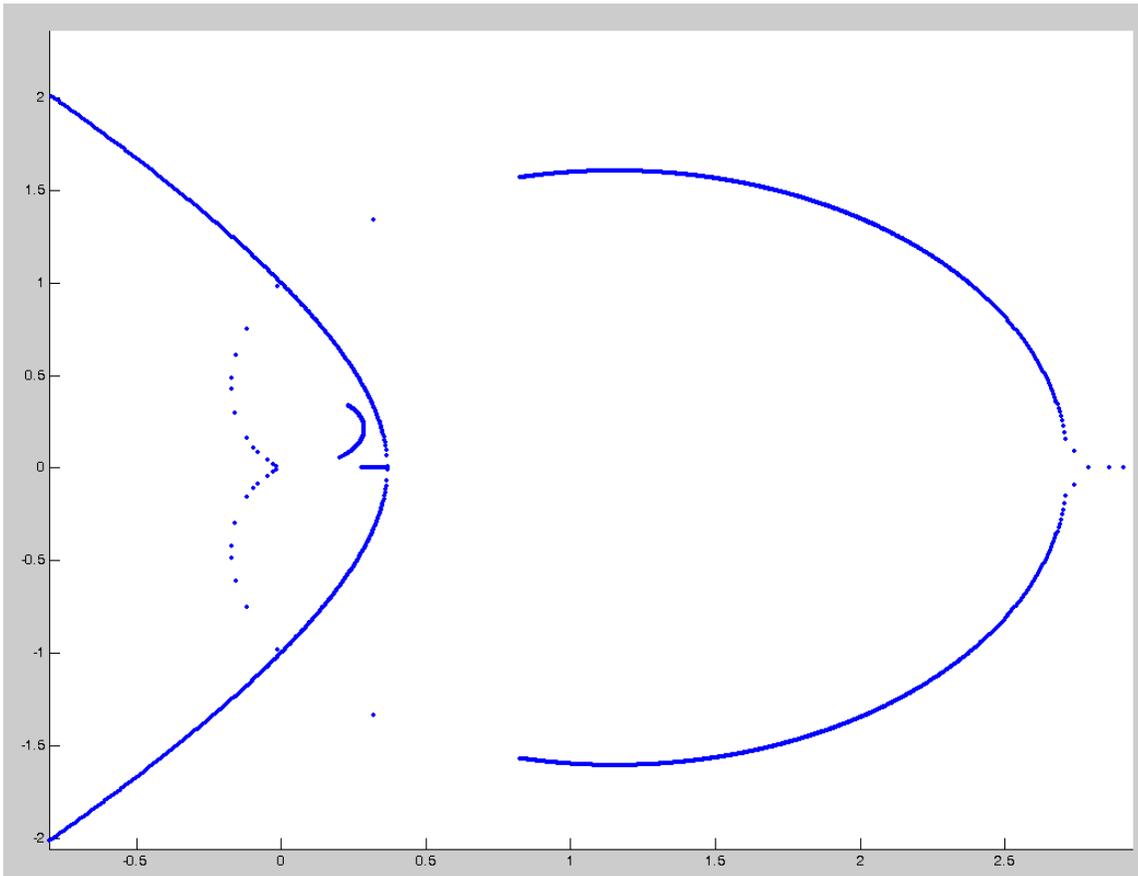

Figure 11

**Discussion and Conclusions**
Given these results, it stands to reason that one can hit any given attractor in the complex plane with the right base (complex, with both a real and an imaginary part), but that these attractors are arranged in ordered structures that follow as a function of the base, sometimes also as a function of the starting domain (positive or negative i's), whereas the specific starting point does not seem to matter in any of the cases we considered.

It is a matter of debate whether it is unfortunate or fortuitous that the level of investigation in this paper remained on the level of empirical mathematics, as opposed to an analytical approach that deduced these structural properties from the axiomatic characteristics of the logarithmic operation itself. This also implies that the meaning of these structures is not easily assessed.

One obvious problem is that these results – showing singularity-like attractor points – could be owed to the fact that all numerical methods, including the ones in this paper suffer from a numerical representation that is necessarily finite in terms of its precision. Thus, numerical methods are inherently vulnerable to rounding errors. However, while this is possible in principle, it is unlikely to account for these results. First, they appear extremely systematic and ordered; second, repeating all calculations with variable precision had no impact on the results.

Finally, we would like to note that it is curious that the attractor seems to be iteratively approached in the shape of a logarithmic spiral (e.g. Figure 1 or 2). Moreover, it is very curious that the liftoff into the imaginary parts happens at an attractor value of about e (e.g. Table 3). In addition, it is interesting that we were able to elicit purely mathematical singularities with an interesting geometrical structure that do not arise from physics. Nor do they arise from complexity theory or nonlinear dynamics, which is usually the case when singularity and attractor systems do arise in pure mathematics. Also, it should be noted that the entire universe of real numbers can be mapped onto these well defined non-trivial (non-zero) singularity points via simple iterative logs.

Hence, there might be possible applications to this method – for instance in cryptography or information compression – as one aspect of the information (the starting point) is completely lost by the iterative process, whereas the information about the base is completely preserved, due to the structural properties of these singularities in the complex plane.